\documentclass[12pt,a4paper]{article}
\usepackage{amsmath,amssymb, amsthm}
\newtheorem{theorem}{Theorem}[section]
\newtheorem{lemma}[theorem]{Lemma}
\newtheorem{corollary}[theorem]{Corollary}
\newtheorem{proposition}[theorem]{Proposition}
\numberwithin{equation}{section} 
\makeatletter
\renewcommand\section{\@startsection{section}{1}{\z@}%
                                  {-3.5ex \@plus -1ex \@minus -.2ex}%
                                  {2.3ex \@plus.2ex}%
                                  {\normalfont\bfseries}}
\makeatother

\renewcommand{\leq}{\leqslant}
\renewcommand{\geq}{\geqslant}

\newcommand{\Rkg}{R_{KG}}
\newcommand{\Rkc}{R_{KC}}
\newcommand{\Rkh}{R_{KH}}
\newcommand{\Rkf}{R_{KF}}
\newcommand{\Rkgt}{R_{K\tilde{G}}}
\newcommand{\Rkct}{R_{K\tilde{C}}}

\newcommand{\ind}{\stackrel{\rm \scriptscriptstyle ind}{=}}
\newcommand{\proj}{\stackrel{\rm \scriptscriptstyle proj}{=}}
\newcommand{\rad}{{\rm rad}}
\newcommand{\downa}{\!\!\downarrow}
\newcommand{\upa}{\!\!\uparrow}
\begin{document}
\parindent=0.3in
\baselineskip=18pt plus 1pt
\title{Periodicity of Adams operations on the Green ring of a finite group}
\author{R. M. Bryant and Marianne Johnson}
\maketitle
\bibliographystyle{amsplain}
\begin{abstract}
The Adams operations $\psi_\Lambda^n$ and $\psi_S^n$ on the Green
ring of a group $G$ over a field $K$ provide a framework for the
study of the exterior powers and symmetric powers of $KG$-modules.
When $G$ is finite and $K$ has prime characteristic $p$ we show that
$\psi_\Lambda^n$ and $\psi_S^n$ are periodic in $n$ if and only if
the Sylow $p$-subgroups of $G$ are cyclic. In the case where $G$ is
a cyclic $p$-group we find the minimum periods and use recent work
of Symonds to express $\psi_S^n$ in terms of $\psi_\Lambda^n$.
\end{abstract}
\section{Introduction}
In this paper we study the Adams operations defined from the
exterior powers $\Lambda^n(V)$ and the symmetric powers $S^n(V)$ of
a $KG$-module $V$, where $G$ is a group and $K$ is a field. By a
$KG$-module we shall always mean a finite-dimensional right
$KG$-module.
The modules $\Lambda^n(V)$ and $S^n(V)$ are of fundamental
importance in the study of $KG$-modules. Thus, for a given $V$,
we would like to know the structure of $\Lambda^n(V)$
and $S^n(V)$. For example, we would like to be able to express
these modules up to isomorphism as direct sums of indecomposables.

When describing modules up to isomorphism it is useful to work in
the Green ring (or representation ring) $\Rkg$. This has a
$\mathbb{Z}$-basis consisting of representatives of the isomorphism
classes of indecomposable $KG$-modules, and multiplication comes
from tensor products. Each $KG$-module may be regarded, up to
isomorphism, as an element of $\Rkg$.

For each $n>0$ there are $\mathbb{Z}$-linear maps $\psi_\Lambda^n$
and $\psi_S^n$ from $\Rkg$ to $\Rkg$, called the Adams operations
(see Section 2). The module $\Lambda^n(V)$ may be expressed, in
$\mathbb{Q}\otimes_{\mathbb{Z}}\Rkg$, as a polynomial in
$\psi_\Lambda^1(V), \ldots, \psi_\Lambda^n(V)$. Similarly, $S^n(V)$
is a polynomial in $\psi_S^1(V), \ldots, \psi_S^n(V)$. Thus
$\Lambda^n(V)$ and $S^n(V)$ are determined, up to isomorphism, by
knowledge of the Adams operations.

The advantage of working with the Adams operations is that they have
some rather simple properties (see Section 2). For example, when $n$
is not divisible by the characteristic of $K$, we have
$\psi_\Lambda^n=\psi_S^n$ and $\psi_\Lambda^n$ is a ring endomorphism of
$\Rkg$. Furthermore, according to the evidence available in \cite{BJ}
and in some of the references cited there,
expressions for $\psi_\Lambda^n(V)$ and $\psi_S^n(V)$ within $\Rkg$
are often much simpler in form than expressions for $\Lambda^n(V)$
and $S^n(V)$.

In a previous paper~\cite{BJ} we studied the Adams operations for a
cyclic $p$-group $C$ in prime characteristic $p$. We showed how to
calculate $\psi_\Lambda^n$ for all $n$ not divisible by $p$ (when we
also have $\psi_\Lambda^n=\psi_S^n$). Here we shall continue our
study of $\psi_\Lambda^n$ and $\psi_S^n$ on $\Rkc$, allowing $n$ to
be an arbitrary positive integer. In Sections 4 and 6 we
shall establish periodicity in
$n$: explicitly, $\psi_\Lambda^n=\psi_\Lambda^{n+2q}$ and
$\psi_S^n=\psi_S^{n+2q}$, where $q=|C|$. The proofs of these results
make use of our earlier work~\cite{BJ}, and the proof for $\psi_S^n$
also relies upon deep work of Symonds~\cite{Sy} that gives
periodicity `modulo induced modules' of the symmetric powers of
indecomposable $KC$-modules. Symonds~\cite{Sy} also obtained a
result that expresses symmetric powers of indecomposable
$KC$-modules in terms of exterior powers, but again `modulo induced
modules'. In Section 6 we show
that there is a corresponding result that
expresses $\psi_S^n$ directly in terms of $\psi_\Lambda^n$.

The other main results of this paper concern the periodicity of
$\psi_\Lambda^n$ and $\psi_S^n$ for an arbitrary finite group $G$
in prime characteristic $p$.
We show in Section 3
that $\psi_\Lambda^n$ is periodic in $n$ if and only if the
Sylow $p$-subgroups of $G$ are cyclic. In Section 5 we
obtain a similar
result for $\psi_S^n$, but the proof is much harder and makes use
of~\cite{Sy}. In Section 7 we obtain some lower bounds for
the minimum periods.

\section{Preliminaries}
\noindent Let $G$ be a group and $K$ a field. We consider
(finite-dimensional right)
$KG$-modules and denote the associated Green ring by $\Rkg$,
as in Section 1.
For any $KG$-module $V$, we also write $V$ for the corresponding
element of $\Rkg$. Thus, for $KG$-modules $V$ and $W$, we have $V=W$
in $\Rkg$ if and only if $V \cong W$. The elements $V + W$ and $VW$
of $\Rkg$ correspond to the modules $V \oplus W$ and $V \otimes_K W$,
respectively. The identity element of $\Rkg$ is the one-dimensional
$KG$-module on which $G$ acts trivially, usually denoted by $1$ when
regarded as an element of $\Rkg$ and $K$ when regarded as a
$KG$-module.

If $V$ is a $KG$-module and $n$ is a non-negative integer, then we
regard the $n$th exterior power $\Lambda^n(V)$ and the $n$th
symmetric power $S^n(V)$ as elements of $\Rkg$. In particular,
$\Lambda^0(V)=S^0(V)=1$ and $\Lambda^1(V)=S^1(V)=V$. If  $V$ has
dimension $r$ then $\dim \Lambda^n(V) =
\binom{r}{n}$, where for $n>r$ we define $\binom{r}{n}=0$. Thus we
see that $\Lambda^r(V)$ is a one-dimensional module, whilst
$\Lambda^n(V)=0$ for all $n>r$.

For any $KG$-module $V$ we define elements of the power-series ring
$\Rkg[[t]]$ by
$$\Lambda(V,t) = 1 + \Lambda^1(V)t +\Lambda^2(V)t^2 + \cdots,$$
$$S(V,t) = 1 + S^1(V)t +S^2(V)t^2 + \cdots .$$
(In fact $\Lambda(V,t)$ belongs to the polynomial ring
$\Rkg[t]$.) We extend $\Rkg$ to a ring $\mathbb{Q}\Rkg$,
allowing coefficients from $\mathbb{Q}$, so that $\mathbb{Q}\Rkg \cong
\mathbb{Q}\otimes_{\mathbb{Z}} \Rkg$. The formal expansion of
$\log(1+x)$ yields elements $\log\Lambda(V,t)$ and $\log S(V,t)$
of $\mathbb{Q}\Rkg[[t]]$. Thus for $n>0$ we may define elements
$\psi_\Lambda^n(V)$ and $\psi_S^n(V)$ of $\mathbb{Q}\Rkg$ by the
equations
\begin{equation}
\psi_\Lambda^1(V) - \psi_\Lambda^2(V)t + \psi_\Lambda^3(V)t^2 -
\cdots = \displaystyle \frac{d}{dt}
\log\Lambda(V,t) = \Lambda(V,t)^{-1} \displaystyle \frac{d}{dt}
\Lambda(V,t),
\label{AdamsdefLambda}
\end{equation}
\begin{equation}
\psi_S^1(V) + \psi_S^2(V)t + \psi_S^3(V)t^2 + \cdots =
\displaystyle\frac{d}{dt}\log S(V,t) = S(V,t)^{-1} \displaystyle
\frac{d}{dt} S(V,t).
\label{AdamsdefS}
\end{equation}

\vspace{0.1in}
\noindent
It is easily verified that $\psi_\Lambda^1(V) = \psi_S^1(V) =V$ for
all $KG$-modules $V$. In fact, it can be shown that
$\psi_\Lambda^n(V), \psi_S^n(V) \in \Rkg$ for all $n>0$. For the
trivial module $K$ we have $\psi_\Lambda^n(K) = \psi_S^n(K) = K$ for
all $n > 0$. Furthermore,
$$\psi_\Lambda^n(V+W) = \psi_\Lambda^n(V) + \psi_\Lambda^n(W),\;\;\;
\psi_S^n(V+W) = \psi_S^n(V) + \psi_S^n(W),$$
for all $n >0$ and all $KG$-modules $V$ and $W$ (see~\cite{Br} for
details). Thus the definitions of $\psi_\Lambda^n$ and $\psi_S^n$
may be extended to give $\mathbb{Z}$-linear maps
$$\psi_\Lambda^n : \Rkg \rightarrow \Rkg,\;\;\; \psi_S^n :
\Rkg \rightarrow \Rkg,$$
called the $n$th Adams operations on $\Rkg$. We  write $\delta$ for
the `dimension' map $\delta : \Rkg \rightarrow \mathbb{Z}$. This is
the ring homomorphism satisfying $\delta(V) = \dim V$ for
every $KG$-module $V$. Since there is an embedding
$\nu: \mathbb{Z} \to \Rkg$ given by $\nu(k) = k.1$ for all
$k \in \mathbb{Z}$ we may also regard $\delta$ as an endomorphism
of $\Rkg$. By \cite[(2.2)]{BJ} we have
\begin{equation}
\label{dimemsion}
\delta(\psi_\Lambda^n(A))=\delta(\psi_S^n(A))=\delta(A),
\end{equation}
for all $A\in \Rkg$ and all $n>0$.

We now extend $\Rkg$ to a ring $\mathbb{C}\Rkg$, allowing
coefficients from $\mathbb{C}$, so that $\mathbb{C}\Rkg \cong
\mathbb{C}\otimes_{\mathbb{Z}} \Rkg$ and $\mathbb{C}\Rkg$ is a
$\mathbb{C}$-algebra. The Adams operations on
$\Rkg$ extend to $\mathbb{C}$-linear maps $\psi_\Lambda^n$ and
$\psi_S^n$
on $\mathbb{C}\Rkg$. The
following result may be found in \cite[Theorem 5.4 and (4.4)]{Br}.

\begin{proposition}
\label{factorisation} Let $G$ be a group, $K$ a field, and $n$ a
positive integer not divisible by the characteristic of $K$.
Then $\psi_\Lambda^n = \psi_S^n$ and $\psi_\Lambda^n$ is a
ring endomorphism of $\Rkg$ and an algebra
endomorphism of $\mathbb{C}\Rkg$. If\/ $V$ is a $KG$-module then
$\psi_\Lambda^n(V)$ is a
${\mathbb{Z}}$-linear combination of direct summands of $V^{\otimes n}$.
Furthermore, under composition of maps we have
$\psi_\Lambda^n \circ \psi_\Lambda^{n'} = \psi_\Lambda^{nn'}$ and
$\psi_S^n \circ \psi_S^{n'} = \psi_S^{nn'}$ for every positive integer
$n'$.
\end{proposition}

For any $KG$-module $V$,
\eqref{AdamsdefLambda} yields
$$\biggl( \sum_{n=1}^{\infty}(-1)^{n-1}\psi_\Lambda^n(V)t^{n-1}\biggr)
\biggl( \sum_{n=0}^{\infty}\Lambda^n(V)t^n \biggr) =
\sum_{n=1}^\infty n\Lambda^n(V)t^{n-1}.$$
By comparing coefficients of $t^{n-1}$ we obtain Newton's formula
\begin{eqnarray}
\sum_{j=0}^{n-1} (-1)^{j+1} \psi_\Lambda^{n-j}(V)\Lambda^j(V)&=&
(-1)^{n}n\Lambda^n(V) \label{NewtonLambda1}
\end{eqnarray}
for all $n>0$ and all $KG$-modules $V$. Hence
\begin{eqnarray}
\psi_\Lambda^n(V) &=& \sum_{j=1}^{n-1}(-1)^{j+1}
\psi_\Lambda^{n-j}(V)\Lambda^j(V)\; + (-1)^{n+1}n\Lambda^n(V).
\label{NewtonLambda2}
\end{eqnarray}
A similar equation may be obtained for $\psi_S^n(V)$ by
means of \eqref{AdamsdefS}.

\begin{lemma}
\label{Adamsinflate} Let $K$ be a field of prime characteristic $p$,
$G$ a finite $p$-group, $V$ the regular
$KG$-module and $n$ a positive integer not divisible by $p$. Then
$$\psi_\Lambda^n(V) = \psi_S^n(V) = V.$$
\end{lemma}

\noindent
{\bf Proof.}
Since $G$ is a $p$-group, $V$
is indecomposable and $V^{\otimes n}$ is a direct multiple of $V$.
Therefore, by Proposition
\ref{factorisation} and \eqref{dimemsion}, we have
$\psi_\Lambda^n(V) = \psi_S^n(V) = V$. \hfill $\square$

\vspace{0.2in}
From now on we assume that $G$ is a finite group. For any subgroup
$H$ of $G$ and any $KG$-module $V$ we write $V\downa_H$ to
denote the $KH$-module obtained from $V$ by restriction. Since
restriction commutes with taking direct sums and tensor products
there is a ring homomorphism from $\Rkg$ to
$\Rkh$ mapping $V$ to $V\downa_H$ for every
$KG$-module $V$. For each $A \in \Rkg$ we write $A\downa_H$ for
the image of $A$ under this homomorphism. Since restriction
commutes with
the formation of exterior and symmetric powers, we have
(see \cite[Lemma 2.3]{Br2})
\begin{equation}
\label{Adamsrestrict} \psi_\Lambda^n(A)\downa_H \ =
\psi_\Lambda^n(A\downa_H), \;\;\;\psi_S^n(A)\downa_H \ =
\psi_S^n(A\downa_H),
\end{equation}
for all $n>0$ and all $A \in \Rkg$. (Note that the Adams
operations on the left-hand sides of these equations act on
$\Rkg$ whilst those on the right act on $\Rkh$.)

For any $KH$-module $U$ we write $U\upa^G$ to denote the
$KG$-module obtained from $U$ by induction. Since induction commutes
with taking direct sums there is a $\mathbb{Z}$-linear
map from $\Rkh$ to $\Rkg$ that maps $U$ to $U\upa^G$ for every
$KH$-module $U$. For each $A \in \Rkh$ we write $A\upa^G$ for the
image of $A$ under this map.
For any $KG$-module $V$ and any $KH$-module $U$ we have
$U\upa^GV = (U (V\downa_H))\upa^G$ in $\Rkg$ (by~\cite[Proposition
3.3.3]{Be}, for example). Hence the $\mathbb{Z}$-span of all modules
induced from $H$ is an ideal of $\Rkg$. Similarly, the
$\mathbb{Z}$-span of all relatively $H$-projective $KG$-modules is
an ideal of $\Rkg$.
In particular, the $\mathbb{Z}$-span
of all projective modules is an ideal of $\Rkg$. For $A,B \in \Rkg$
we write $A \proj B$ to mean that $A - B$ is  a
$\mathbb{Z}$-linear combination of projectives.

For a $KG$-module $V$ let $\mathcal{P}(V)$ denote a projective cover
of $V$ (see~\cite[Section 1.5]{Be} for further details).
The Heller translate
$\Omega(V)$ is then defined (up to isomorphism) as the kernel of the
map $\mathcal{P}(V) \twoheadrightarrow V$. The definition of
$\Omega$ extends to give a
$\mathbb{Z}$-linear map $\Omega: \Rkg \rightarrow \Rkg$.
Hence, for each $n \geq 0$, we have a $\mathbb{Z}$-linear map
$\Omega^n: \Rkg \rightarrow
\Rkg$ defined by composition: thus $\Omega^0$ is the
identity map and $\Omega^n = \Omega \circ \Omega^{n-1}$ for $n > 0$.
It can be shown that $\Omega(AB) \proj
\Omega(A)B$ for all $A, B\in \Rkg$ (see~\cite[Corollary 3.1.6]{Be}).
Repeated application of this
result yields
\begin{equation}
\label{Hellermult} \Omega^i(A)\Omega^j(B) \proj \Omega^{i+j}(AB)
\end{equation}
for all $i,j \geq 0$ and all $A, B\in \Rkg$.

For any subgroup $H$ of $G$ the map $\Omega: \Rkh \to \Rkh$ is defined
in the same way as $\Omega: \Rkg \to \Rkg$ (but using $KH$-modules
instead of $KG$-modules). It follows easily from Schanuel's lemma
(see \cite[Lemma 1.5.3]{Be}) that $\Omega(V)\downa_H \ \proj
\Omega(V\downa_H)$ for every $KG$-module $V$. Thus, for all
$A \in \Rkg$, we have
\begin{equation}
\label{Schanuel} \Omega(A)\downa_H \ \proj \Omega(A\downa_H).
\end{equation}

Let $N$ be a normal subgroup of $G$. Each $K(G/N)$-module yields a
$KG$-module by `inflation'. Indeed, inflation yields a ring
embedding $\mu: R_{K(G/N)} \rightarrow \Rkg$. Also, since $\mu$
commutes with the formation of exterior and symmetric powers, we have
\begin{equation}
\label{inflation} \psi_\Lambda^n(\mu(A))=\mu(\psi_\Lambda^n(A)),
\;\;\;\psi_S^n(\mu(A))=\mu(\psi_S^n(A)),
\end{equation}
for all $n > 0$ and all $A \in R_{K(G/N)}$ (see \cite[Lemma 2.3]{Br2}).

\section{Periodicity of the Adams operations $\psi_\Lambda^n$}
\noindent
Let $G$ be a finite group and $K$ a field. One of the main themes of
this paper is the periodicity, as a function of $n$, of the Adams
operations $\psi_\Lambda^n$ and $\psi_S^n$ on $\Rkg$. When we
refer to the periodicity of $\psi_\Lambda^n$ or $\psi_S^n$ we shall
always mean periodicity in $n$.

\begin{lemma}
\label{charzero} Let $G$ be a finite group and let $K$ be a field
such that $|G|$ is
not divisible by the characteristic of\/ $K$.
Then $\psi_\Lambda^n =
\psi_\Lambda^{n+e}$ and $\psi_S^n = \psi_S^{n+e}$, for all $n>0$,
where $e$ is the exponent of\/ $G$.
\end{lemma}

\noindent
{\bf Proof.}
(When ${\rm char}\,K =0$ this result is essentially well known.) For any
$KG$-module $V$ let ${\rm Br}(V)$ denote the Brauer character of $V$.
Thus ${\rm Br}(V): G \to \mathbb{C}$. Hence we define
${\rm Br}(A): G \to \mathbb{C}$ for all $A \in \Rkg$ so that
${\rm Br}$ is $\mathbb{Z}$-linear on $\Rkg$.
Let $V$ be a $KG$-module. Then, by \cite[Lemma 2.6]{Br2}, we have
$${\rm Br}(\psi_\Lambda^n(V))(g) = {\rm Br}(V)(g^n) =
{\rm Br}(\psi_S^n(V))(g)$$
for all $g \in G$. Hence
$${\rm Br}(\psi_\Lambda^n(V)) = {\rm Br}(\psi_\Lambda^{n+e}(V)) =
{\rm Br}(\psi_S^{n+e}(V))={\rm Br}(\psi_S^n(V)).$$
However, for all $A, B \in \Rkg$ we have ${\rm
Br}(A) ={\rm Br}(B)$ if and only if $A=B$ in $\Rkg$ (as can be derived
from \cite[Corollary 5.3.6]{Be}). Hence the result follows.
\hfill $\square$

\vspace{0.2in}
Under the assumptions of Lemma \ref{charzero}, $e$ is, in fact, the
minimum period of $\psi_\Lambda^n$ and $\psi_S^n$, as we shall see
in Section 7. When ${\rm char}\,K = 0$, Lemma \ref{charzero} yields the
periodicity of $\psi_\Lambda^n$ and $\psi_S^n$
for every finite group. Thus for the rest of this
section we shall assume that $K$ has prime characteristic $p$. Here we
concentrate on $\psi_\Lambda^n$. The periodicity of $\psi_S^n$
will be studied in Section 5.

It was proved in \cite[Theorem 7.2]{Br2} that $\psi_\Lambda^n$ and
$\psi_S^n$ are periodic (with specified upper bounds for the periods) when
the Sylow $p$-subgroups of $G$ have order $p$. Here we generalise
this fact for $\psi_\Lambda^n$ as follows. The
proof is relatively simple
but gives only a crude upper bound for the period.

\begin{theorem}
\label{mainext} Let $G$ be a finite group and let $K$ be a field of
prime characteristic $p$. Then the Adams operations $\psi_\Lambda^n$
on the Green ring $\Rkg$ are periodic in $n$ if and only if the
Sylow $p$-subgroups of\/ $G$ are cyclic.
\end{theorem}

\noindent
{\bf Proof.}
Suppose first that the Sylow $p$-subgroups of $G$ are cyclic. By a
theorem of Higman (see~\cite[Theorem (64.1)]{CR1}), there
are only finitely many isomorphism classes of indecomposable
$KG$-modules. Hence $\mathbb{C}\Rkg$ is finite-dimensional. By work
of Green and O'Reilly (see~\cite[Theorem (81.90)]{CR2}),
$\mathbb{C}\Rkg$ is semisimple. Thus $\mathbb{C}\Rkg$ has a
$\mathbb{C}$-basis $\{e_i: i=1, \ldots, m\}$ consisting of
pairwise-orthogonal idempotents.
Note that every idempotent of $\mathbb{C}\Rkg$ has the form
$\sum_{j \in S} e_j$, where $S$ is a subset of $\{1,\dots,m\}$, and
every endomorphism of $\mathbb{C}\Rkg$ maps each $e_i$ to an
idempotent. Thus $\mathbb{C}\Rkg$ has only finitely many
endomorphisms.
By Proposition \ref{factorisation},
$\psi_\Lambda^n$ is an endomorphism of $\mathbb{C}\Rkg$ for all $n$
not divisible by $p$. Thus there are only finitely many
possibilities for the maps $\psi_\Lambda^n$ where $p \nmid n$.

Choose a positive integer $d$ such that $\dim V \leq p^d$ for every
indecomposable $KG$-module $V$. If $n$ is a positive integer not
divisible by $p^{d+1}$ we may write $n=p^ik$ where $0 \leq i \leq d$
and $p \nmid k$. Thus, by Proposition \ref{factorisation},
$\psi_\Lambda^n = \psi_\Lambda^k \circ \psi_\Lambda^{p^i}$. There
are only finitely many possibilities for $\psi_\Lambda^k$ and at most
$d+1$ possibilities for $\psi_\Lambda^{p^i}$. Hence there are only
finitely many possibilities for the maps $\psi_\Lambda^n$ where $p^{d+1}
\nmid n$.

If $c=p^dk$ where $p\nmid k$ then $p^{d+1} \nmid n$ for all $n \in
\{c, c+1, \ldots, c+p^d-1\}$. Hence there are only finitely many
possibilities for the $p^d$-tuple $(\psi_\Lambda^c, \ldots,
\psi_\Lambda^{c+p^d-1})$ where $c$ has the given form. It follows
that there exist positive integers $a$ and $b$ such that $a<b$ and
\begin{equation}
\label{repeatuple} (\psi_\Lambda^a, \ldots, \psi_\Lambda^{a+p^d-1})
= (\psi_\Lambda^b, \ldots, \psi_\Lambda^{b+p^d-1}).
\end{equation}
Write $s=b-a$. We shall show that $\psi_\Lambda^n =
\psi_\Lambda^{n+s}$ for all $n>0$.

Let $V$ be any indecomposable $KG$-module and set $r=\dim V$. Thus
$r \leq p^d$. It suffices to show that $\psi_\Lambda^n(V) =
\psi_\Lambda^{n+s}(V)$ for all $n>0$. For each $n$ define an
$r$-tuple $\Psi_n$ by
$$\Psi_n = (\psi_\Lambda^n(V), \ldots, \psi_\Lambda^{n+r-1}(V)). $$
It suffices to show that $\Psi_n = \Psi_{n+s}$ for all $n>0$. By
\eqref{repeatuple}, we have $\Psi_a = \Psi_{a+s}$. Thus it suffices
to show, for all $n>0$, that $\Psi_n = \Psi_{n+s}$ if and only if
$\Psi_{n+1} = \Psi_{n+1+s}$.

Since $\Lambda^j(V)=0$ for $j>r$, Newton's formulae,
\eqref{NewtonLambda1} and \eqref{NewtonLambda2}, with $n$ replaced
by $n+r$, become
\begin{equation}
\label{Newton(n+r)1} \sum_{j=0}^r (-1)^{j+1} \psi_\Lambda^{n+r-j}(V)
\Lambda^j(V) = 0
\end{equation}
and
\begin{equation}
\label{Newton(n+r)2} \psi_\Lambda^{n+r}(V) = \sum_{j=1}^r (-1)^{j+1}
\psi_\Lambda^{n+r-j}(V) \Lambda^j(V),
\end{equation}
for all $n>0$. Since $\Lambda^r(V)$ is one-dimensional,
there exists a one-dimensional $KG$-module $W$
such that $\Lambda^r(V)W=1$ in $\Rkg$. Hence, from
\eqref{Newton(n+r)1}, we obtain
\begin{equation}
\label{Newton(n+r)3} (-1)^r\psi_\Lambda^{n}(V) = \sum_{j=0}^{r-1}
(-1)^{j+1} \psi_\Lambda^{n+r-j}(V) \Lambda^j(V)W.
\end{equation}

The elements $\psi_\Lambda^{n+r-j}(V)$ on the right of
\eqref{Newton(n+r)2} are the components of $\Psi_n$ (in reverse
order). Hence if $\Psi_n=\Psi_{n+s}$ we obtain
$\psi_\Lambda^{n+r}(V)=\psi_\Lambda^{n+r+s}(V)$ and hence
$\Psi_{n+1}=\Psi_{n+1+s}$. Similarly the elements
$\psi_\Lambda^{n+r-j}(V)$ on
the right of \eqref{Newton(n+r)3} are the components of
$\Psi_{n+1}$. Hence if $\Psi_{n+1}=\Psi_{n+1+s}$ we obtain
$\psi_\Lambda^{n}(V)=\psi_\Lambda^{n+s}(V)$ and hence
$\Psi_n=\Psi_{n+s}$. This gives the required result.

We shall now prove the converse, omitting some of the details for
the sake of brevity.
Let $G$ be a finite group with
non-cyclic Sylow $p$-subgroups. We prove that $\psi_\Lambda^n$
is non-periodic. In fact we prove the stronger result that the
$\psi_\Lambda^n$ are distinct for $p \nmid n$. (Since
$\psi_\Lambda^n=\psi_S^n$ when $p \nmid n$, our proof shows that
$\psi_S^n$ is also non-periodic.)

Let $H$ be a minimal non-cyclic $p$-subgroup of $G$. Thus either $H
\cong C_p \times C_p$ or $p=2$ and $H \cong Q_8$, where $C_p$
has order $p$ and $Q_8$ is the quaternion group.

Suppose first that $H \cong C_p \times C_p$. Write $K$ for the
trivial one-dimensional $KH$-module, the identity element of $\Rkh$.
The Heller translates $\Omega^n(K)$, for $n \geq 1$, are isomorphic
to the kernels of the maps in a minimal projective resolution of
$K$. Since $H \cong C_p \times C_p$, this resolution is
non-periodic. (If it is periodic its modules have bounded dimension
and so the spaces ${\rm Ext}_{KH}^n(K,K)$ have bounded dimension. But
${\rm Ext}_{KH}^n(K,K) \cong H^n(H,K)$, and $H^n(H,K)$ has dimension $n+1$
by the K{\"u}nneth formula.) It follows that the $\Omega^n(K)$ are
distinct elements of $\Rkh$ for $n \geq1$. Furthermore, since $K$ is
a non-projective indecomposable, each $\Omega^n(K)$ is a
non-projective indecomposable. Write $V=\Omega(K)$. Then, by
\eqref{Hellermult}, we have
$$\Omega^n(K) = \Omega^n(K^n) \proj (\Omega(K))^n = V^n,$$
for all $n \geq 1$. Thus $V^{\otimes n}$ is isomorphic to a direct sum of
$\Omega^n(K)$ and projectives.

Let $n$ be a positive integer not divisible by $p$. Then, by
Proposition \ref{factorisation}, $\psi_\Lambda^n(V)$ is a
$\mathbb{Z}$-linear combination of direct summands of $V^{\otimes n}$.
However, $\dim V = p^2 - 1$ and so
$\delta(\psi_\Lambda^n(V)) = p^2-1$
by \eqref{dimemsion}. Thus $\psi_\Lambda^n(V)$ is not a
$\mathbb{Z}$-linear combination of projectives only and
$\psi_\Lambda^n(V)$ must involve $\Omega^n(K)$. Therefore, for
values of $n$ not
divisible by $p$, the elements $\psi_\Lambda^n(V)$ of $\Rkh$ are
distinct.

Let $U$ be the Heller translate of the trivial one-dimensional
$KG$-module. By \eqref{Schanuel} we have
$U\downa_H \ \proj V = \Omega(K)$. Since
$\Omega(K)$ is a non-projective indecomposable, we
find that $U\downa_H \ =V+W$, where $W$ is projective.

Let $m$ and $n$ be distinct positive integers not divisible by $p$. Then, as
seen above, $\psi_\Lambda^m(V) \neq \psi_\Lambda^n(V)$. However, by
Lemma \ref{Adamsinflate}, $\psi_\Lambda^m(W) = \psi_\Lambda^n(W)$.
Therefore $\psi_\Lambda^m(U\downa_H) \neq
\psi_\Lambda^n(U\downa_H)$. Hence $\psi_\Lambda^m(U) \neq
\psi_\Lambda^n(U)$, by \eqref{Adamsrestrict}. Thus $\psi_\Lambda^m
\neq \psi_\Lambda^n$ on $\Rkg$.

Now suppose that $p=2$ and $H\cong Q_8$. Thus $H$ has a normal
subgroup $N$ such that $H/N \cong C_2 \times C_2$. Let $V$ be the
$K(H/N)$-module defined by $V=\Omega(K)$, as above.
Thus $V$ is three-dimensional and is isomorphic to the augmentation
ideal of $K(H/N)$.
As we have seen, the
elements $\psi_\Lambda^n(V)$ of $R_{K(H/N)}$ are distinct for $p
\nmid n$. We now regard $V$ as a $KH$-module by inflation. Thus
$N$ acts trivially on $V$ and, by
\eqref{inflation}, the elements $\psi_\Lambda^n(V)$ of $\Rkh$ are
distinct for $p \nmid n$.

Let $M=V\upa^G$. By Mackey's decomposition theorem
(see~\cite[Theorem 3.3.4]{Be}), $M\downa_H$ is isomorphic to a
direct sum of modules $V(g)$, where $g$ ranges over a set of
representatives of double cosets $HgH$. Indeed,
$V(g) = (V\otimes g)\downa_{H^g \cap H}\;\upa^H$,
where $V\otimes g$ is the $K(H^g)$-module `conjugate' to the
$KH$-module $V$ in which we have
$(v \otimes g)g^{-1}hg=vh \otimes g$
for all $v \in V$, $h \in H$.

Note that the subgroups of $H$ are $H$, $\{1\}$, $N$, and cyclic
subgroups $L_1, L_2, L_3$ of order $4$, all of these subgroups being
normal. Note also that $H^g$ has only one element of order $2$
and this acts trivially on $V \otimes g$.
It is straightforward to check that if $H^g \cap H = H$ then
$V(g) \cong V$. If $H^g \cap H =\{1\}$ then $V(g)$ is clearly a
free $KH$-module. If $H^g \cap H = N$ then $N$ acts trivially on
$V\otimes g$ and so $V(g)$ is a free $K(H/N)$-module (regarded as a
$KH$-module). If $H^g \cap H = L_i$ for $i \in \{1,2,3\}$ then
it is straightforward to check that $(V\otimes
g)\downa_{H^g \cap H}$ is the direct sum of a one-dimensional trivial
$KL_i$-module and a regular $K(L_i/N)$-module: hence $V(g)$ is the
direct sum of a regular $K(H/L_i)$-module and a regular
$K(H/N)$-module.
It follows that $M\downa_H \ = rV +W$, where $r$ is a positive
integer and $W$ is a sum of regular modules for factor groups
of $H$ regarded as $KH$-modules.

Let $m$ and $n$ be distinct positive integers not divisible by $p$. Then
we have $\psi_\Lambda^m(rV) \neq \psi_\Lambda^n(rV)$. Also, by Lemma
\ref{Adamsinflate} and \eqref{inflation}, we have
$\psi_\Lambda^m(W) = \psi_\Lambda^n(W)$. It follows that
$\psi_\Lambda^m(M) \neq \psi_\Lambda^n(M)$, and so
$\psi_\Lambda^m \neq \psi_\Lambda^n$ on $\Rkg$.
\hfill $\square$

\section{The Adams operations $\psi_\Lambda^n$ for
a cyclic $p$-group}
\noindent In this section we shall consider the Adams operations
$\psi_\Lambda^n$ on the Green ring $\Rkc$, where $K$ is a field of
prime characteristic $p$ and $C$ is a cyclic $p$-group of order $q
\geq1$. As we have seen in Theorem \ref{mainext},
$\psi_\Lambda^n$ is periodic in $n$. Here we shall show, for $q > 1$,
that the minimum period is $2q$. We shall also establish the
symmetry property that $\psi_\Lambda^n = \psi_\Lambda^{2q-n}$
for all $n < 2q$.

There are precisely $q$ indecomposable $KC$-modules up to isomorphism,
which we denote by $V_1, \ldots, V_q$, where $V_r$ has dimension $r$
for $r=1, \ldots, q$. Thus $\Rkc$ has $\mathbb{Z}$-basis
$\{V_1,\dots,V_q\}$.
The one-dimensional module $V_1$ is
the identity element of $\Rkc$, which we shall sometimes write simply
as $1$. When $q=1$ we have $\Rkc=\mathbb{Z}V_1$ so that
$\psi_\Lambda^n$ is the identity map for all $n > 0$. Thus we
usually assume that $q > 1$.

For $q>1$ let $\tilde{C}$ denote the subgroup of index $p$ in $C$.
Thus $\tilde{C}$ is cyclic of order $q/p$. We write $\tilde{V}_1,
\ldots, \tilde{V}_{q/p}$ to denote the indecomposable
$K\tilde{C}$-modules, up to isomorphism, where $\tilde{V}_r$ has
dimension $r$, for $r=1, \ldots, q/p$. Let $r \in \{1,\dots,q\}$.
Then it is well known (and easy to prove) that
\begin{equation}
\label{restriction} V_r\downa_{\tilde{C}} \ = (p-b)\tilde{V}_a + b
\tilde{V}_{a+1},
\end{equation}
where $r=ap+b$ with $0 \leq b <p$. (We take the convention that
$\tilde{V}_0 = 0$.) Notice that $a+1 \leq q/p$, provided that $b
>0$. Also, for $r = 1,\dots,q/p$, we have
\begin{equation}
\label{induction} \tilde{V}_r\upa^C \ = V_{pr}.
\end{equation}
A $KC$-module will be said to be induced if it is induced from a
$K\tilde{C}$-module, and an element of $\Rkc$ will be said to be
induced if it is a $\mathbb{Z}$-linear combination of induced
modules. The set of induced elements is
an ideal of $\Rkc$
and, by \eqref{induction}, this ideal has $\mathbb{Z}$-basis
$\{V_p, V_{2p}, \ldots, V_q\}$. For $A,B
\in \Rkc$ we write $A \ind B$ to mean that $A - B$ is induced.

\begin{lemma}
\label{resind} Let $q>1$ and $A \in \Rkc$.

\noindent{\rm (i)} If\/ $A \ind 0$ and
$A\downa_{\tilde{C}} \ = 0$ then $A=0$.

\noindent{\rm (ii)} If\/ $A \ind 0$ and
$A\downa_{\tilde{C}} \  \proj 0$ then $A \proj 0$.
\end{lemma}

\noindent
{\bf Proof.} (i) Suppose that $A \ind 0$ and
$A\downa_{\tilde{C}} \ = 0$.
Since $A \ind 0$ we may write
$A = \sum_{r=1}^{q/p} \alpha_r V_{pr}$,
where $\alpha_r \in \mathbb{Z}$ for $r=1, \ldots, q/p$. Since
$A\downa_{\tilde{C}} \ = 0$ we have
$\sum_{r=1}^{q/p} \alpha_r p\tilde{V}_{r} = 0$ by \eqref{restriction}.
Thus $\alpha_r = 0$
for all $r$, and hence $A=0$.

(ii) This is similar.
\hfill $\square$

\vspace{0.2in}
The regular module $V_q$ is the unique projective indecomposable
$KC$-module. Thus, for $A,B \in \Rkc$, we have $A \proj B$ if and
only if  $A - B \in \mathbb{Z}V_q$. For $q>1$ we have that $V_q$
is induced and so $A \proj B$ implies $A \ind B$ for all
$A,B \in \Rkc$.

If $1 \leq p^j \leq q$ then $C$ has a factor group $C(p^j)$ of order
$p^j$ and inflation of modules gives a ring embedding
$R_{KC(p^j)} \rightarrow \Rkc$. Thus, for $r = 1,\dots,p^j$, we may
identify $V_r$ with the indecomposable $KC(p^j)$-module of
dimension $r$. Since $V_{p^j}$ is the only projective indecomposable
$KC(p^j)$-module we have
\begin{equation}
\label{Pmult} \mbox{$V_r V_{p^j} = rV_{p^j}$, \, for $r = 1,\dots,p^j$.}
\end{equation}

Let $P_{KC}$ denote the $\mathbb{Z}$-span in $\Rkc$ of the set of
all permutation modules for $C$ over $K$. Each transitive
permutation module $M$ is the module induced from the one-dimensional
trivial module for some subgroup of $C$ (namely, the stabilizer of a
point). Thus, by repeated use of \eqref{induction} (applied to subgroups
of $C$ rather than $C$ itself), we find that $M \cong V_{p^j}$ for some $j$
such that $1 \leq p^j \leq q$. Conversely, the modules $V_{p^j}$
are permutation modules. Thus $P_{KC}$ has $\mathbb{Z}$-basis
$\{V_1, V_p, V_{p^2},\dots,V_q\}$. Since the tensor product of permutation
modules is a permutation module, $P_{KC}$ is a subring of $\Rkc$.
(This also follows from \eqref{Pmult}.)

For $A \in \Rkc$ write $A = \sum_{i=1}^{q}\alpha_i(A)V_i$, where
$\alpha_i(A) \in \mathbb{Z}$ for $i=1, \ldots, q$. Thus $\alpha_i(A)$
denotes the multiplicity with which $V_i$ occurs in $A$.
The following result is an easy consequence of Lemma
\ref{resind}\,(i).

\begin{corollary}
\label{resalpha} Let $q>1$ and let $A, B \in P_{KC}$.
If\/ $\alpha_1(A) = \alpha_1(B)$ and
$A\downa_{\tilde{C}} \ =B\downa_{\tilde{C}}$ then $A=B$.
\end{corollary}

Let $r \in \{1,\dots,q\}$ and $j \in \{0,\dots,r\}$. Then
$\dim \Lambda^j(V_r) = \dim \Lambda^{r-j}(V_r)$. In fact
\begin{equation}
\label{palindromic} \Lambda^j(V_r) = \Lambda^{r-j}(V_r)
\end{equation}
in $R_{KC}$. This is well known but we sketch a proof.
Since $\dim
\Lambda^r(V_r)=1$ we have $\Lambda^r(V_r) \cong K$ as a $KC$-module. Also
there is an isomorphism of
$K$-spaces, $\theta: \Lambda^j(V_r) \rightarrow {\rm Hom}_K
(\Lambda^{r-j}(V_r), \Lambda^r(V_r))$, induced by multiplication in the
exterior algebra $\Lambda(V_r)$. We may regard ${\rm Hom}_K
(\Lambda^{r-j}(V_r), \Lambda^r(V_r))$ as a $KC$-module, the
contragredient dual of $\Lambda^{r-j}(V_r)$. Then it is easily
verified that $\theta$ is an isomorphism of $KC$-modules. Hence
\eqref{palindromic} follows because all
$KC$-modules are self-dual.

\begin{lemma}
\label{Lambda(reg)} For each non-negative integer $n$ we have
$\Lambda^n(V_q) \in P_{KC}$ and
\begin{equation*}
\alpha_1(\Lambda^n(V_q)) = \begin{cases}
1& \mbox{if\, $n=0$ or $n=q$,}\\
0&\mbox{otherwise.}
\end{cases}
\end{equation*}
\end{lemma}

\noindent
{\bf Proof.}
Since $\Lambda^0(V_q)=\Lambda^q(V_q)=V_1$ and $\Lambda^n(V_q)=0$ for
$n>q$, the result holds for $n=0$ and for all $n \geq q$. Thus we
assume that $1 \leq n < q$.

Let $g$ be a generator of $C$. Then the regular module $V_q$ has
basis $\{x_1, \ldots, x_q\}$ where
$x_ig = x_{i+1}$ for $1 \leq i < q$ and $x_qg=x_1$. It is easily
checked that the set
$$\{(-1)^{(i_1+\cdots + i_n)(n+1)} x_{i_1} \wedge x_{i_2}
\wedge \cdots \wedge x_{i_n}: 1 \leq i_1 < i_2 < \cdots < i_n \leq q\}$$
is a basis of $\Lambda^n(V_q)$ invariant under the action of $g$.
Thus $\Lambda^n(V_q)$ is a permutation module, and so
$\Lambda^n(V_q) \in P_{KC}$. Since $n < q$ there are no orbits of
length one in the given basis. Hence $\alpha_1(\Lambda^n(V_q)) = 0$.
\hfill $\square$

\begin{lemma}
\label{psi(reg)1} Let $q>1$ and $n > 0$.

\noindent {\rm (i)} {\it We have\/
$\psi_\Lambda^n(V_q) \in P_{KC}$.}

\noindent {\rm(ii)} {\it If\/ $p$ is odd then }
$$\alpha_1(\psi_\Lambda^n(V_q)) = \begin{cases}
q& \mbox{if\, $q \mid n$,}\\
0&\mbox{otherwise. }
\end{cases}$$

\noindent {\rm(iii)} {\it If\/ $p=2$ then}
$$\alpha_1(\psi_\Lambda^n(V_q)) = \begin{cases}
(-1)^{n/q}q & \mbox{if\, $q \mid n$,}\\
0&\mbox{otherwise. }
\end{cases}$$
\end{lemma}

\noindent
{\bf Proof.}
We prove the result by induction on $n$.
The result is trivial for $n=1$ because $\psi_\Lambda^1(V_q) = V_q$.
Now let $n>1$ and suppose that the result holds for all $k$ such
that $1 \leq k <n$. By \eqref{NewtonLambda2} we have that
\begin{equation}
\label{Newton(reg)} \psi_\Lambda^n(V_q) = \sum_{j=1}^{n-1}
(-1)^{j+1} \psi_\Lambda^{n-j}(V_q) \Lambda^j(V_q)
+(-1)^{n+1}n\Lambda^n(V_q).
\end{equation}
By our inductive hypothesis, $\psi_\Lambda^{n-j}(V_q) \in P_{KC}$ for
$j=1, \ldots, n-1$. Also, by Lemma \ref{Lambda(reg)},
$\Lambda^j(V_q) \in P_{KC}$ for $j=1, \ldots, n$. Since $P_{KC}$ is a subring
of $\Rkc$ it now follows from \eqref{Newton(reg)} that
$\psi_\Lambda^n(V_q) \in P_{KC}$, as required for (i).

For all $k>0$ let $\beta_k = \alpha_1(\Lambda^k(V_q))$ and
$\gamma_k = \alpha_1(\psi_\Lambda^k(V_q))$. Then, by
\eqref{Newton(reg)} and \eqref{Pmult}, we have
$$\gamma_n = \sum_{j=1}^{n-1} (-1)^{j+1}
\gamma_{n-j} \beta_j +(-1)^{n+1} n \beta_n.$$
By Lemma \ref{Lambda(reg)} it follows that
$\gamma_n = 0$ if $n < q$, $\gamma_q = (-1)^{q+1}q$ and
$\gamma_n = (-1)^{q+1}\gamma_{n-q}$ if $n > q$.
Parts (ii) and (iii) now follow from our inductive hypothesis.
\hfill $\square$

\vspace{0.2in}
For all positive integers $a$ and $b$, let $(a,b)$ denote the
greatest common divisor of $a$ and $b$.

\begin{proposition}
\label{psi(reg)2} Let $n>0$.

\noindent {\rm(i)} {\it If\/ $p$ is odd then\/ $\psi_\Lambda^n(V_q)
= (n,q)V_{q/(n,q)}.$}

\noindent {\rm(ii)} {\it If\/ $p=2$ then}
$$\psi_\Lambda^n(V_q) = \begin{cases}
V_q& \mbox{if\, $n$ is odd,}\\
(n,2q)V_{2q/(n,2q)} - (n,q)V_{q/(n,q)} &\mbox{if\, $n$ is even. }\\
\end{cases}$$
\end{proposition}

\noindent
{\bf Proof.}
When $q=1$, (i) and (ii) hold trivially. From now on we assume that
$q>1$ and write $\tilde{q}=q/p$. Define $A_n = (n,q)V_{q/(n,q)}$
and $\tilde{A}_n =
(n,\tilde{q})\tilde{V}_{\tilde{q}/(n,\tilde{q})}$.
By considering the cases where $q \mid n$ and $q \nmid n$
separately, it is easily verified that
\begin{equation}
\label{A(n,q)restrict} A_n\downa_{\tilde{C}} \ =
p\tilde{A}_n.
\end{equation}
It is also easy to see that
\begin{equation}
\label{A(n,q)alpha} \alpha_1(A_n) = \begin{cases}
q& \mbox{if\, $q \mid n$,}\\
0&\mbox{otherwise.}
\end{cases}
\end{equation}

(i) Suppose that $p$ is odd. It suffices to show that
$\psi_\Lambda^n(V_q)=A_n$. By Lemma \ref{psi(reg)1}\,(ii) and
\eqref{A(n,q)alpha}  we have
$\alpha_1(\psi_\Lambda^n(V_q)) = \alpha_1(A_n)$.
By \eqref{Adamsrestrict} and \eqref{restriction} we have
$\psi_\Lambda^n(V_q)\downa_{\tilde{C}} \ =
p\psi_\Lambda^n(\tilde{V}_{\tilde{q}})$. Also we may assume
by induction on $q$ that
$\psi_\Lambda^n(\tilde{V}_{\tilde{q}}) = \tilde{A}_n$. Thus, by
\eqref{A(n,q)restrict}, we have
$\psi_\Lambda^n(V_q)\downa_{\tilde{C}} \ = A_n\downa_{\tilde{C}}$.
However, $\psi_\Lambda^n(V_q) \in P_{KC}$, by Lemma \ref{psi(reg)1}\,(i),
and, clearly, $A_n \in P_{KC}$.
The result now follows by Corollary \ref{resalpha}.

(ii) Suppose that $p=2$. Define
$B_n = V_q$ if $n$ is odd and $B_n = 2A_{n/2} - A_n$ if $n$ is even,
and define $\tilde{B}_n$ similarly.
Then $B_n\downa_{\tilde{C}} \ = 2\tilde{B}_n$, by
\eqref{A(n,q)restrict}, and it suffices to show that
$\psi_\Lambda^n(V_q)=B_n$.
Separating the cases $(n,2q) = 2q$, $(n,2q) = q$ and
$(n,2q) < q$ in Lemma \ref{psi(reg)1}\,(iii) and
\eqref{A(n,q)alpha}, we obtain
$\alpha_1(\psi_\Lambda^n(V_q)) = \alpha_1(B_n)$.
By \eqref{Adamsrestrict} and \eqref{restriction} we have
$\psi_\Lambda^n(V_q)\downa_{\tilde{C}} \ =
2 \psi_\Lambda^n(\tilde{V}_{\tilde{q}})$. Also we may assume
by induction that
$\psi_\Lambda^n(\tilde{V}_{\tilde{q}}) = \tilde{B}_n$. Thus
$\psi_\Lambda^n(V_q)\downa_{\tilde{C}} \ = B_n\downa_{\tilde{C}}$.
The result follows by Corollary \ref{resalpha}.
\hfill $\square$

\vspace{0.2in}
In Sections 6 and 7 we shall also need the corresponding result for
$\psi_S^n(V_q)$.

\vspace{0.1in}
\begin{proposition}
\label{psi(reg)3} For all $n>0$ we have $\psi_S^n(V_q) =
(n,q)V_{q/(n,q)}$.
\end{proposition}

\noindent
{\bf Proof.}
It is easy to verify that $S^n(V_q)$ is a permutation module with
respect to the usual basis. Counting orbits of length one yields that
$\alpha_1(S^n(V_q)) = 1$ if $q \mid n$ and
$\alpha_1(S_n(V_q)) = 0$ otherwise.
By arguments similar to the proofs of Lemma \ref{psi(reg)1}\,(i) and (ii),
we deduce that $\psi_S^n(V_q) \in P_{KC}$ where
$\alpha_1(\psi_S^n(V_q)) = q$ if $q \mid n$ and
$\alpha_1(\psi_S^n(V_q)) = 0$ otherwise.
The result then follows by an argument similar to
the proof of Proposition \ref{psi(reg)2}.
\hfill $\square$

\begin{proposition}
\label{psi^notqs} For all $n>0$ such that $q \nmid n$ we have
$\psi_\Lambda^n = \psi_\Lambda^{n+2p(n,q)}$.
\end{proposition}

\noindent
{\bf Proof.}
If $(n,q)=1$ then,
by \cite[Theorem 3.3]{BJ}, we have $\psi_\Lambda^n =
\psi_\Lambda^{n+2p}$, as required.
So suppose that $1 < (n,q) < q$. Then $n=(n,q)k$, where $p \nmid k$, and
hence $n+2p(n,q) = (n,q)(k+2p)$, where $p\nmid k+2p$. Thus
$\psi_\Lambda^n = \psi_\Lambda^k \circ \psi_\Lambda^{(n,q)}$ and
$\psi_\Lambda^{n+2p(n,q)} = \psi_\Lambda^{k+2p} \circ
\psi_\Lambda^{(n,q)}$, by Proposition \ref{factorisation}.
Since $p\nmid k$ we have $\psi_\Lambda^k = \psi_\Lambda^{k+2p}$.
Thus $\psi_\Lambda^n = \psi_\Lambda^{n + 2p(n,q)}$.
\hfill $\square$

\vspace{0.2in}
Recall that $\delta$ is the endomorphism of $\Rkc$
satisfying $\delta(V_r) = rV_1$ for all $r$.

\begin{theorem}
\label{periodL} Let $K$ be a field of prime characteristic
$p$ and let $C$ be a cyclic $p$-group of order $q>1$.

\noindent {\rm(i)} The Adams operations $\psi_\Lambda^n$ on $\Rkc$
are periodic in $n$, with minimum period $2q$.

\noindent {\rm(ii)} We have $\psi_\Lambda^n = \psi_\Lambda^{2q-n}$
for all $n<2q$.

\noindent {\rm(iii)} We have
$\psi_\Lambda^{2q} = \delta. $
\end{theorem}

\noindent
{\bf Proof.}
(i) We begin by proving, by induction on $n$, that
$\psi_\Lambda^n = \psi_\Lambda^{n+2q}$
for all $n>0$. By Proposition \ref{psi^notqs},
the result holds for all $n$ such that $q \nmid n$.
In particular the result holds for all $n<q$. Now
assume that $n \geq q$ and that $\psi_\Lambda^k =
\psi_\Lambda^{k+2q}$ for all $k<n$. Let $r \in \{1,\ldots, q\}$. It is
enough to prove that $\psi_\Lambda^n(V_r) =
\psi_\Lambda^{n+2q}(V_r)$.

For $r=q$ this follows easily from Proposition \ref{psi(reg)2}.
Thus we may take $r<q$. Since $n \geq q > r$ and $\Lambda^j(V_r) = 0$
for $j > r$,
\eqref{NewtonLambda2} gives
$$
\psi_\Lambda^n(V_r) = \sum_{j=1}^{r} (-1)^{j+1}
\psi_\Lambda^{n-j}(V_r)\Lambda^j(V_r)
$$
and
$$
\psi_\Lambda^{n+2q}(V_r) = \sum_{j=1}^{r} (-1)^{j+1}
\psi_\Lambda^{n+2q-j}(V_r)\Lambda^j(V_r).
$$
Thus, by
our inductive hypothesis, $\psi_\Lambda^n(V_r) =
\psi_\Lambda^{n+2q}(V_r)$.

We have shown that $\psi_\Lambda^n$ is
periodic in $n$, with minimum period $\lambda$ dividing $2q$. It remains
to prove that $\lambda =2q$. First suppose that $p=2$. By Proposition
\ref{psi(reg)2}\,(ii), we see that $\psi_\Lambda^q(V_q) =
\psi_\Lambda^{q+s}(V_q)$ if and only if $s$ is a multiple of $2q$.
Hence $2q$ divides $\lambda$, giving $\lambda=2q$.
Now suppose that $p$ is odd.
By Proposition \ref{psi(reg)2}\,(i), we see that $\psi_\Lambda^q(V_q)
= \psi_\Lambda^{q+s}(V_q)$ if and only if $s$ is a multiple of $q$.
Hence $q$ divides $\lambda$, giving that $\lambda =q$ or
$\lambda =2q$. Since $p$ is
odd, $q+1$ is even. Thus, by~\cite[Theorem 5.1]{BJ}, we
find that $\psi_\Lambda^1(V_2) \neq \psi_\Lambda^{q+1}(V_2)$. Hence
$\lambda = 2q$.

(ii) The result is trivial for $n=q$. By symmetry, it
remains only to prove the result for $n<q$. Then we may write
$n=(n,q)k$ where $(n,q)<q$ and $p \nmid k$. Also, $2q-n = (n,q)k'$,
where $k'=2q/(n,q) - k$. Since $p$ divides $2q/(n,q)$ we have $p
\nmid k'$. Thus by Proposition \ref{factorisation} we have
$\psi_\Lambda^n = \psi_\Lambda^k \circ \psi_\Lambda^{(n,q)}$
and $\psi_\Lambda^{2q-n} = \psi_\Lambda^{k'} \circ
\psi_\Lambda^{(n,q)}$.
Since $k' \equiv -k \;(\bmod \; 2p)$, we have
$\psi_\Lambda^k = \psi_\Lambda^{k'}$ by~\cite[Corollary 3.4]{BJ}.
Thus $\psi_\Lambda^n=\psi_\Lambda^{2q-n}$.

(iii) It suffices to prove
that $\psi_\Lambda^{2q}(V_r) = rV_1$ for $r \in \{1,\dots,q\}$.
By \eqref{NewtonLambda2},
$$\psi_\Lambda^{2q}(V_r) = \sum_{j=1}^r(-1)^{j+1}
\psi_\Lambda^{2q-j}(V_r)\Lambda^j(V_r),$$
and, by (ii), we have
$\psi_\Lambda^{2q-j} = \psi_\Lambda^j$ for all $j$ such that $1 \leq j
< 2q$. Moreover, by \eqref{palindromic}, we have $\Lambda^j(V_r)
= \Lambda^{r-j}(V_r)$ for all $j$ such that $1 \leq j \leq r$. Hence
$$\psi_\Lambda^{2q}(V_r) = \sum_{j=1}^r(-1)^{j+1}
\psi_\Lambda^j(V_r)\Lambda^{r-j}(V_r)=(-1)^r\sum_{k=0}^{r-1}(-1)^{k+1}
\psi_\Lambda^{r-k}(V_r)\Lambda^k(V_r). $$
Thus, by \eqref{NewtonLambda1}, we see that $\psi_\Lambda^{2q}(V_r) =
r \Lambda^r(V_r) = rV_1$, as required.
\hfill $\square$

\section{Periodicity of the Adams operations $\psi_S^n$}
\noindent In this section we investigate the periodicity of $\psi_S^n$
for the Green ring of an arbitrary finite
group $G$ over a field $K$. By Lemma 3.1, $\psi_S^n$ is periodic
if the characteristic of $K$ does not divide $|G|$.
Thus we assume that $K$ has prime
characteristic $p$. We shall show that $\psi_S^n$ is periodic
if and only if $G$ has cyclic Sylow $p$-subgroups.

The main step in the proof is to establish that $\psi_S^n$
is periodic in the special
case when $G$ has a normal cyclic Sylow $p$-subgroup with cyclic
factor group. We assume that $G$ has this form for almost all of
this section: only in Theorem \ref{mainsym} will we go to the
general situation, using Conlon's induction theorem. We now fix some
hypotheses and notation that will apply until Theorem \ref{mainsym}.

Assume that $K$ is algebraically closed.
Let $G$ be a finite group with a normal cyclic Sylow $p$-subgroup
$C$ such that $G/C$ is cyclic. Thus (by the Schur-Zassenhaus
theorem) $G$ has a cyclic $p'$-subgroup $H$ such that $G=HC$. Write
$q=|C|$ and $m=|G/C|=|H|$, and let $C=\langle g \rangle$ and $H=\langle h
\rangle$. By Theorem \ref{mainext} there exists a positive integer
$\pi$ such that $\psi_\Lambda^n = \psi_\Lambda^{n+\pi}$ holds on
$R_{KF}$ for all $n>0$ and for every subgroup $F$ of $G$.
However,
when $p \nmid n$, we have $\psi_\Lambda^n=\psi_S^n$. Thus we have
the following result.

\begin{lemma}
\label{properties} There exists a positive integer $\pi$ satisfying

\noindent
{\rm (i)} $\pi$ is divisible by $p$ and $m$,\\
\noindent {\rm (ii)} $\psi_S^n=\psi_S^{n+\pi}$ on $R_{KF}$
for all $n$ such
that $p \nmid n$ and for every subgroup $F$ of $G$.
\end{lemma}

All $KH$-modules will be regarded as
$KG$-modules via the epimorphism $G \twoheadrightarrow H$ that
restricts to the identity on $H$. Thus $\Rkh$ becomes a subring
of $\Rkg$. We consider the indecomposable $KG$-modules, as
described in \cite[pp.\ 34--37, 42--43]{Al}.
We summarise the main facts, but adapt the presentation
in~\cite{Al} and use different notation. (In particular we use right
modules instead of left modules.)

First note that, since $H$ is cyclic and $K$ is algebraically
closed, there is a faithful one-dimensional $KH$-module $X$ and the
irreducible $KH$-modules up to isomorphism are the tensor powers $1,
X, \ldots, X^{m-1}$. (We shall usually use Green ring notation for
modules.) Clearly $X$ has multiplicative order $m$ in $\Rkh$. The
irreducible $KG$-modules are the same as the irreducible
$KH$-modules.

When $q>1$ let $\kappa$ be an integer such that $h^{-1}gh=g^\kappa$
and reduce $\kappa$ modulo $p$ to obtain a non-zero element
$\bar{\kappa}$ of the base field of $K$. Let $W$ be the
one-dimensional $KH$-module on which $h$ acts as the scalar
$\bar{\kappa}$. (This coincides with $W$ as used in
\cite{Al} and as described in
\cite[Exercise 5.3]{Al}.) Also, let $d$ be the
multiplicative order of $\bar{\kappa}$. Thus $d \mid m$ and
$d \mid p-1$. By
suitable choice of $X$ we can take $W=X^{m/d}$. Note that $W$ has
multiplicative order $d$ in $\Rkh$. When $q=1$ we take $d=1$
and $W=1$.

Let $J_q$ be the projective cover of the
one-dimensional trivial $KG$-module.
Then (see~\cite[p.\ 37]{Al}) $J_q$ is uniserial with proper
submodules $J_q(g-1)^r$, for $r=1, \ldots, q$, and composition
factors $1, W, W^2, \ldots, W^{q-1}$, from top to bottom. For $r=1,
\ldots, q$, write $J_r = J_q/J_q(g-1)^r$. It is easily verified that
$J_r\downa_C \ = V_r$ in the notation of Section 4. For $i=0,
\ldots, m-1$, the (tensor product) module $X^iJ_r$ is uniserial, and
its composition factors from top to bottom are $X^i, X^iW, \ldots,
X^iW^{r-1}$. The modules $X^iJ_r$, for $i=0, \ldots, m-1$ and $r=1,
\ldots, q$, are the indecomposable $KG$-modules up to isomorphism
(see~\cite[p.\ 42]{Al}). Note also that $(X^iJ_r)\downa_C \ = V_r$.

When $q>1$ let $\tilde{C}$ denote the subgroup of index $p$ in $C$,
that is, $\tilde{C} = \langle g^p \rangle$. For $s=1, \ldots, q/p$,
let $\tilde{V}_s$ denote the indecomposable $K\tilde{C}$-module
of dimension $s$, as in Section 4.

\begin{lemma}
\label{relproj} When $q>1$ the indecomposable $KG$-modules that are
projective relative to $\tilde{C}$ are those of the form $X^iJ_{ps}$
with $0 \leq i \leq m-1$ and $1 \leq s \leq q/p$.
\end{lemma}

\noindent
{\bf Proof.}
The module $X^iJ_{ps}$ is relatively $C$-projective
(by~\cite[Theorem 9.2]{Al}, for example). Hence $X^iJ_{ps}$ is a
summand of $V_r\upa^G$ for some $r$. By Mackey's decomposition
theorem, $V_r\upa^G \;\downa_C \ = mV_r$, and so we must have that
$r=ps$. However, $V_{ps} = \tilde{V}_s\upa^C$. Therefore $X^iJ_{ps}$
is a summand of $\tilde{V}_s\upa^G$ and is relatively
$\tilde{C}$-projective.

Conversely, suppose that $X^iJ_r$ is relatively
$\tilde{C}$-projective. Then $X^iJ_r$ is a summand of
$\tilde{V}_s\upa^G$ for some $s$. Hence $X^iJ_r$ is a summand of
$V_{ps}\upa^G$. Thus, by Mackey's decomposition theorem, $r=ps$.
\hfill $\square$

\vspace{0.2in}
When $q>1$ let $\tilde{G}=H\tilde{C} < G$ and let $\tilde{J}_1,
\ldots, \tilde{J}_{q/p}$ be the $K\tilde{G}$-modules defined in the
same way as we defined $J_1, \ldots, J_q$ for $KG$. Thus the
indecomposable $K\tilde{G}$-modules are the modules $X^i\tilde{J}_s$
for $i=0, \ldots, m-1$ and $s=1, \ldots, q/p$.

\begin{lemma}
\label{relprojres} When $q>1$ we have that
$J_{ps}\!\downa_{\tilde{G}} \ =
(1+W + \cdots + W^{p-1})\tilde{J}_s$ for
$s=1, \ldots, q/p$.
\end{lemma}

\noindent
{\bf Proof.}
Write $M=J_{ps}\downa_{\tilde{G}}$ and let $X^i\tilde{J}_k$ be any
indecomposable summand of $M$. Since $M\downa_{\tilde{C}} \
= p\tilde{V}_s$, we must have $k=s$. Thus we may write
$M=U\tilde{J}_s$, where $U$ is a $KH$-module. We have $M/\rad(M)=U$
in $\Rkgt$, because $\tilde{J}_s/\rad(\tilde{J}_s)$ is trivial.
However, $\rad(M) = M(g^p-1)$, by \cite[Lemma 5.8]{Al} applied
to $\tilde{G}$. Since $J_{ps}(g^p-1) = J_{ps}(g-1)^p$, the
composition factors of $J_{ps}/J_{ps}(g^p-1)$ are $1, W, \ldots,
W^{p-1}$. However
$(J_{ps}/J_{ps}(g^p-1))\downa_{\tilde{G}} \ = M/\rad(M)$ and
$M/\rad(M)$ is completely reducible. Thus
$U = 1+W+\cdots+W^{p-1}$.
Hence $M=(1+W+ \cdots
+W^{p-1})\tilde{J}_s$, as required.
\hfill $\square$

\vspace{0.2in}
Assume that $q>1$. In the Green ring $\Rkg$ the set $I$ of all
$\mathbb{Z}$-linear combinations of relatively
$\tilde{C}$-projective $KG$-modules is an ideal (see Section 2). For
$A,B \in \Rkg$ we write $A \approx B$ to mean that $A-B \in I$. Note
that when $G=C$ we have $A \approx B$ if and only if
$A \ind B$. Thus the following result is a generalisation of Lemma
\ref{resind}\,(i).

\begin{lemma}
\label{resind2} Let $q>1$ and let $A$ be an element of\/ $\Rkg$ such that
$A \approx 0$ and
$A\downa_{\tilde{G}} \ =0$. Then $A=0$.
\end{lemma}

\noindent
{\bf Proof.}
By Lemma \ref{relproj} we may write $A=\sum_s U_s J_{ps}$ where $s$ ranges
over $\{1, \ldots, q/p\}$ and $U_s \in \Rkh$
for all $s$. By Lemma \ref{relprojres} and the assumption
$A\downa_{\tilde{G}} \ =0$ we obtain
$$\sum_s U_s(1+W + \cdots + W^{p-1})\tilde{J}_s = 0$$
in $\Rkgt$. By the remark before Lemma \ref{relprojres} it follows
that we have
\begin{equation}
\label{zerodivisor} U_s(1+W + \cdots +W^{p-1})=0
\end{equation}
in $\Rkh$ for all $s \in \{1,\dots,q/p\}$.
It suffices to prove that $U_s=0$. Recall that $W$ has
multiplicative order $d$, where $d \mid p-1$. Thus, by
\eqref{zerodivisor}, we have
\begin{eqnarray}
\label{Usrewritten}
U_s &=&-U_s(1+W + \cdots + W^{p-2})\nonumber\\
&=& -((p-1)/d)U_s(1+W + \cdots + W^{d-1}).
\end{eqnarray}
Recall that $\{1, X, \ldots,
X^{m-1}\}$ is a basis of $\Rkh$. Since $W=X^{m/d}$ this basis
may be written as
$\{X^iW^j: 0 \leq i \leq l-1,\, 0 \leq j \leq d-1\}$, where
$l = m/d$. Hence there are integers $n_{i,j}$ such that
$$-((p-1)/d)U_s = \sum_{i=0}^{l-1}\sum_{j=0}^{d-1} n_{i,j}X^iW^j.$$
Since $W^j(1+ W + \cdots +
W^{d-1}) = (1+ W + \cdots + W^{d-1})$ for all $j$, we see, by
\eqref{Usrewritten}, that there exist integers $n_0, n_1, \ldots,
n_{l-1}$ such that
\begin{equation}
\label{form} U_s= \sum_{i=0}^{l-1}n_iX^i(1+W + \cdots + W^{d-1}).
\end{equation}
Let $\phi :\Rkh \rightarrow \Rkh$ be the ring endomorphism defined
by $\phi(X) = X^d$. Then $\phi(W)=X^m=1$. Applying $\phi$ to
\eqref{zerodivisor} yields $p\phi(U_s)=0$ and hence
$\phi(U_s)=0$. Thus, by applying $\phi$ to \eqref{form}, we obtain
the equation $0=d\sum_{i=0}^{l-1}n_iX^{di}$
in $\Rkh$. Therefore $n_i=0$ for all $i$, and so $U_s=0$, as
required.
\hfill $\square$

\vspace{0.2in}
Suppose that $q>1$ and recall that for $A,B \in \Rkg$ we write $A
\approx B$ to mean that $A-B \in I$. Similarly, for $A,B\in
\mathbb{Q}\Rkg$ we write $A \approx B$ to mean that $A-B \in
\mathbb{Q}I$. For $A(t),B(t) \in
\mathbb{Q}\Rkg[[t]]$, where $A(t) = \sum_{i=0}^{\infty}A_it^i$ and
$B(t) = \sum_{i=0}^{\infty}B_it^i$, we write $A(t) \approx B(t)$ to
mean that $A_i \approx B_i$ for all $i$. It is easy to check that if
$A_0=B_0=1$ then $A(t) \approx B(t)$ implies $\log A(t) \approx \log
B(t)$.

Let $V$ be an indecomposable $KG$-module. Thus $V=X^iJ_r$
for some $i$ and $r$. By \cite[Theorem 1.2]{Sy}, there is a
one-dimensional submodule $E$ of $S^q(V)$ such that
$$S(V,t) \approx (1 + B_1t + \cdots + B_{q-1}t^{q-1})(1+Et^q +E^2t^{2q}+
\cdots),$$
where $B_n=S^n(V)$ for $n=1, \ldots, q-1$. Thus, by the remark
above, we have
\begin{equation}
\label{logprod}
\log S(V,t)\approx \log ((1 + B_1t + \cdots + B_{q-1}t^{q-1})
(1+Et^q +E^2t^{2q}+ \cdots)).
\end{equation}
Write $\psi_S(V,t) = \sum_{n = 1}^{\infty} \psi_S^n(V)t^n$. Then,
by \eqref{AdamsdefS}, we have
\begin{eqnarray}
\label{modulo}
\psi_S(V,t) &=& t\displaystyle\frac{d}{dt} \log S(V,t)\nonumber\\
&\approx&  t\displaystyle\frac{d}{dt} \log (1 + B_1t + \cdots +
B_{q-1}t^{q-1}) - t\displaystyle\frac{d}{dt} \log (1-Et^q)\\
&\approx&  \frac{B_1t + 2B_2t^2 + \cdots + (q-1)B_{q-1}t^{q-1}}{1+ B_1t
+ \cdots + B_{q-1}t^{q-1}}+ \frac{qEt^q}{1-Et^q}\nonumber.
\end{eqnarray}
Thus we have
\begin{eqnarray*}
\lefteqn{(\psi_S(V,t)- qEt^q(1-Et^q)^{-1})(1 + B_1t + \cdots +
B_{q-1}t^{q-1})}\\
& & \qquad \qquad \qquad \qquad \qquad \qquad \; \approx B_1t +
2B_2t^2+ \cdots + (q-1)B_{q-1}t^{q-1},
\end{eqnarray*}
so that
$$(D_1t + D_2t^2+\cdots)(1 + B_1t + \cdots + B_{q-1}t^{q-1}) \approx
B_1t + 2B_2t^2+ \cdots + (q-1)B_{q-1}t^{q-1},$$
where
\begin{equation}
\label{Ddef} D_n = \begin{cases}
\psi_S^n(V) & \mbox{if\, $q \nmid n$,}\\
\psi_S^n(V) -qE^{n/q}& \mbox{if\, $q \mid n$.}
\end{cases}
\end{equation}
Thus, for all $n\geq q$, we have
\begin{equation}
\label{BandD} D_n+D_{n-1}B_1 + \cdots + D_{n-q+1}B_{q-1} \approx 0.
\end{equation}

\begin{proposition}
\label{hypoperiodic} Let $K$ be an algebraically closed field
of prime characteristic $p$
and let $G$ be a finite group with a normal cyclic Sylow $p$-subgroup
$C$ such that $G/C$ is cyclic, with $|C|=q$ and $|G/C|=m$. Let $\pi$ be
any positive integer satisfying the conditions
of Lemma\/ {\rm \ref{properties}}.
Then $\psi_S^n=\psi_S^{n+q\pi/p}$ on $\Rkg$ for all $n>0$.
\end{proposition}

\noindent
{\bf Proof.}
First note that $mp \mid \pi$ by Lemma \ref{properties}\,(i).
Thus if $q=1$ the result follows by Lemma \ref{charzero}. Now suppose
that $q>1$. We use the notation introduced above and
we may assume, by induction, that the result holds for
$\tilde{G}=H\tilde{C}<G$. Let $V$ be an indecomposable $KG$-module
and let $n > 0$.
It suffices to show that $\psi_S^n(V) = \psi_S^{n'}(V)$ where
$n'=n+q\pi/p$.

By the inductive hypothesis, $(\psi_S^n(V) -
\psi_S^{n'}(V))\downa_{\tilde{G}} \ = 0$.
We shall show that $\psi_S^n(V) \approx \psi_S^{n'}(V)$.
This gives $\psi_S^n(V) -
\psi_S^{n'}(V) = 0$ by Lemma \ref{resind2}, as required.

We use the notation of \eqref{logprod} and \eqref{Ddef}.
Since $E$ is a one-dimensional $KG$-module, $E$ is a $KH$-module and
$E^m=1$ in $\Rkg$. Since
$\pi/p$ is divisible by $m$ we
have $E^{n/q} = E^{(n+q\pi/p)/q}$ when $q \mid n$.
Hence, by \eqref{Ddef}, the
proposition will follow if we can prove that $D_n \approx D_{n'}$.

First suppose that  $n<q$ and write $n=p^ik$ where $p \nmid k$. Thus
$p^{i+1} \mid q$ and we have
$n' = p^ik'$,
where $k'=k+q\pi/p^{i+1}$. Since $p \mid \pi$
we have $p \nmid k'$. Thus,
by Proposition \ref{factorisation}, $\psi_S^n=\psi_S^k \circ
\psi_S^{p^i}$ and $\psi_S^{n'}=\psi_S^{k'} \circ
\psi_S^{p^i}$. By Lemma \ref{properties}\,(ii), we have
$\psi_S^k=\psi_S^{k'}$. Hence $\psi_S^n(V)=\psi_S^{n'}(V)$,
giving $D_n \approx D_{n'}$ by \eqref{Ddef}.

Now suppose that $n \geq q$. By \eqref{BandD}, we have that $D_n
\approx - \sum_{i=1}^{q-1} D_{n-i}B_i $ and $D_{n'} \approx
-\sum_{i=1}^{q-1} D_{n'-i}B_i$. By induction on $n$ we may assume that
$D_{n-i} \approx D_{n'-i}$ for $i=1, \ldots, q-1$. Hence $D_n\approx
D_{n'}$, as required.
\hfill $\square$

\begin{theorem}
\label{mainsym} Let $G$ be a finite group and let $K$ be a field of
prime characteristic $p$. Then the Adams operations $\psi_S^n$ on
the Green ring $\Rkg$ are periodic in $n$ if and only if the Sylow
$p$-subgroups of $G$ are cyclic.
\end{theorem}

\noindent
{\bf Proof.}
Suppose first that the Sylow $p$-subgroups of $G$ are cyclic.
In order to prove that $\psi_S^n$ is periodic we may assume that $K$
is algebraically closed, by  \cite[Lemma 2.4]{Br2}.
As in~\cite[p.\ 184]{Be}, we say that a finite group is
$p$-hypo-elementary if it is an extension of a $p$-group by a cyclic
$p'$-group. Note that Proposition \ref{hypoperiodic} shows that
$\psi_S^n$ is periodic in $n$ on $\Rkf$ for every
$p$-hypo-elementary subgroup $F$ of $G$. Thus there exists
a positive integer $s$ such that $\psi_S^n$ has period dividing $s$
on $\Rkf$ for all such $F$. We show that
$\psi_S^n(A)=\psi_S^{n+s}(A)$ for all $A \in \Rkg$. By the
choice of $s$, the element $\psi_S^n(V)-\psi_S^{n+s}(V)$ of $\Rkg$
restricts to $0$ in $\Rkf$ for every $p$-hypo-elementary subgroup
$F$ of $G$. Therefore, by~\cite[Corollary 5.6.9]{Be} (a result
obtained from Conlon's induction theorem), we obtain
$\psi_S^n(V)-\psi_S^{n+s}(V)=0$, as required.

The converse was noted in the proof of Theorem \ref{mainext}.
\hfill $\square$

\section{The Adams operations $\psi_S^n$ for a cyclic $p$-group}
\noindent
In this section we consider the operations
$\psi_S^n$ on the Green ring $\Rkc$, where $K$ is a field of prime
characteristic $p$ and $C$ is a cyclic $p$-group. As before we
write $q = |C|$. When $q=1$ we have $\Rkc = {\mathbb{Z}}V_1$
so that $\psi_S^n$ is the identity map for all $n>0$. Thus from now
on we assume that $q>1$.
We establish results for $\psi_S^n$ analogous to those for
$\psi_\Lambda^n$ obtained in Section 4 and, using work of
Symonds \cite{Sy}, we show that $\psi_S^n$ may be expressed
in terms of $\psi_\Lambda^n$.

If $\bar{K}$ is the algebraic closure of $K$ then there is an
isomorphism $R_{KC} \to R_{\bar{K}C}$ mapping $V_r$ to
$\bar{K} \otimes_K V_r$ for $r = 1,\dots,q$. It is easy to see
that this isomorphism commutes with the Heller maps and the Adams
operations. Thus the results in Section 5 obtained under the
assumption that $K$ is algebraically closed hold, in the case $G = C$,
without the need for algebraic closure. We shall use
these results for arbitrary $K$ without further comment.

Our first result is the analogue of
parts (i) and (ii) of Theorem \ref{periodL}. We shall obtain the
analogue of Theorem \ref{periodL}\,(iii) in Corollary \ref{periodL2}
below.

\begin{theorem}
\label{periodS}Let $K$ be a field of prime characteristic $p$ and
let $C$ be a cyclic $p$-group of order $q > 1$.

\noindent {\rm(i)} The Adams operations $\psi_S^n$
on $\Rkc$ are periodic in $n$, with minimum period $\sigma$, where
$\sigma =2q$
when $p$ is odd and $\sigma =q$ when $p=2$.

\noindent {\rm(ii)} We have $\psi_S^n = \psi_S^{\sigma -n}$ for all
$n<\sigma$.
\end{theorem}

\noindent
{\bf Proof.}
(i) By Proposition \ref{hypoperiodic}, $\psi_S^n$
is periodic in $n$. Let $\sigma$ be the minimum period. Then,
by Proposition \ref{psi(reg)3}, we must have $q \mid \sigma$.

Suppose first that $p$ is odd. By \cite[Theorem 3.3]{BJ}
we have $\psi_S^n=\psi_S^{n+2p}$ on $R_{KF}$ for all $n$
such that $p \nmid n$ and for every subgroup $F$ of $C$. Set
$\pi=2p$. Then $\pi$ satisfies the conditions of Lemma
\ref{properties} and hence, by Proposition \ref{hypoperiodic}, we have
that $\psi_S^n = \psi_S^{n+2q}$ on $\Rkc$ for all $n>0$. Thus
$\sigma = q$ or $\sigma = 2q$. By the proof of Theorem
\ref{periodL}\,(i), we have $\psi_S^1 = \psi_\Lambda^1 \neq
\psi_\Lambda^{q+1} = \psi_S^{q+1}$, and so $\sigma = 2q$.

Now suppose that $p=2$. By \cite[Corollary 3.5]{BJ} we have
$\psi_S^n=\psi_S^{n+2}$ on $R_{KF}$ for all $n$ such
that $2 \nmid n$ and for every subgroup $F$ of $C$. Set
$\pi=2$. Then $\pi$ satisfies the conditions of Lemma
\ref{properties} and hence, by Proposition \ref{hypoperiodic}, we have
that $\psi_S^n = \psi_S^{n+q}$ on $\Rkc$ for all $n>0$, giving
$\sigma = q$.

(ii) Let $p$ be arbitrary. By an argument entirely similar
to the proof of Theorem \ref{periodL}\,(ii), it can be shown that
$\psi_S^n = \psi_S^{2q-n}$ for all $n<2q$. However,
$\psi_S^{2q-n} = \psi_S^{\sigma -n}$ for all $n < \sigma$, by (i).
This gives the result.
\hfill $\square$

\vspace{0.2in}
We use all the notation of Section 4 and write
$\tilde{q} = q/p$. Thus $\tilde{q}$ is the order of $\tilde{C}$.
Recall that $\Omega(V)$ denotes the Heller translate of a
$KC$-module $V$ and that $\Omega$ extends to a $\mathbb{Z}$-linear
map $\Omega: \Rkc \rightarrow \Rkc$. It is easy to check that
$\Omega(V_r) = V_{q-r}$, for $r=1,\dots,q$,
where $V_0=0$. Hence $\Omega^2(V_r) \proj V_r$ for $r=1, \ldots,
q$. Thus, working modulo projectives, we see that $\Omega^n$ is
determined by the parity of $n$, with $\Omega^n(V_r) \proj
\Omega(V_r) \proj V_{q-r}$ if $n$ is odd and $\Omega^n(V_r) \proj
V_r$ if $n$ is even.
Similarly we have $\Omega: \Rkct \rightarrow \Rkct$ with
$\Omega(\tilde{V}_r) =
\tilde{V}_{\tilde{q} -r}$, for $r=1,\dots,\tilde{q}$.

Let $r \in \{1,\dots,q\}$ and write $r = ap + b$ where $0 \leq b < p$.
For $b>0$ we have $q-r=p(\tilde{q}-(a+1))
+(p-b)$, where
$0<p-b<p$, while for $b=0$ we have
$q-r=p(\tilde{q}-a)$. Thus, by \eqref{restriction},
\begin{equation}
\label{restrictheller} V_{q-r}\downa_{\tilde{C}} \ =
(p-b)\tilde{V}_{\tilde{q}-a} + b\tilde{V}_{\tilde{q}-(a+1)}.
\end{equation}

Recall that, for $A,B \in \Rkc$, we have $A \proj B$ if and only
if $A-B \in \mathbb{Z}V_q$ while $A \ind B$ if and only if
$A-B \in \mathbb{Z}\{V_p,V_{2p},\dots,V_q\}$. We extend the notation
to $\mathbb{Q}\Rkc$ and $\mathbb{Q}\Rkc[[t]]$ as follows.
For $A,B \in \mathbb{Q}\Rkc$ we write $A \proj B$ when
$A-B \in \mathbb{Q}V_q$ and $A \ind B$ when
$A-B \in \mathbb{Q}\{V_p,V_{2p},\dots,V_q\}$. For
$A(t),B(t) \in \mathbb{Q}\Rkc[[t]]$, where $A(t) =
\sum_{i=0}^\infty A_it^i$ and $B(t) = \sum_{i=0}^\infty B_it^i$,
we write $A(t) \proj B(t)$ when $A_i \proj B_i$ for all $i$ and
$A(t) \ind B(t)$ when $A_i \ind B_i$ for all $i$. Note that
$A(t) \proj B(t)$ implies that $A(t) \ind B(t)$.

We extend the definition of $\Omega$ on $\Rkc$ to a
$\mathbb{Q}$-linear map $\Omega: \mathbb{Q}\Rkc
\rightarrow\mathbb{Q}\Rkc$. Let $J$ be the subring of
$\mathbb{Q}\Rkc[[t]]$ consisting of all elements
$\sum_{i=1}^\infty A_it^i$ with zero constant term and
let $\Omega^*: J \rightarrow J$ be defined by
$$\Omega^*\left(\sum_{i=1}^\infty A_it^i\right) =
\sum_{i=1}^\infty \Omega^i(A_i)t^i$$
for all $\sum_{i=1}^\infty A_it^i \in J$. Clearly
$\Omega^*$ is $\mathbb{Q}$-linear and it is easy to
check, by \eqref{Hellermult}, that $\Omega^*(XY) \proj
\Omega^*(X)\Omega^*(Y)$ for all $X, Y \in J$.
Thus $\Omega^*(XY) \ind \Omega^*(X)\Omega^*(Y)$ for all $X,Y \in J$.
It follows easily that, for all $X \in J$, we have
\begin{equation}
\label{Hellerlog} \Omega^*(\log(1+X)) \ind \log(1+ \Omega^*(X)).
\end{equation}

We shall take $V = V_r$ in \eqref{modulo}, where $r \in \{1,\dots,q\}$.
Note that,
in the notation of \eqref{modulo}, we must have $E = 1$ in $\Rkc$
because $E$ is one-dimensional. Also,
for all $A(t), B(t) \in \mathbb{Q}\Rkc[[t]]$, we have
$A(t) \approx B(t)$ if and only if $A(t) \ind B(t)$,
by the remark preceding Lemma \ref{resind2}. Thus \eqref{modulo} becomes
$$\psi_S(V_r,t) \ind  t\displaystyle\frac{d}{dt} \log (1 + S^1(V_r)t
+ \cdots
+ S^{q-1}(V_r)t^{q-1}) - t\displaystyle\frac{d}{dt} \log (1-t^q).$$

Now let $r \in \{\tilde{q}, \ldots, q\}$. By
\cite[Corollary 3.11]{Sy} we have $S^n(V_r) \ind
\Omega^n(\Lambda^n(V_{q-r}))$ for all $n <q$. Thus
$$\psi_S(V_r,t) \ind  t\displaystyle\frac{d}{dt} \log
\biggl(1 + \sum_{n=1}^{q-1}\Omega^n(\Lambda^n(V_{q-r}))t^n\biggr)
+ q\sum_{k=1}^\infty t^{kq}.$$
Since $\Lambda^n(V_{q-r})=0$ for $n>q-1$, we obtain, by the definition
of $\Omega^*$,
$$\psi_S(V_r,t) \ind
t\displaystyle\frac{d}{dt} \log \biggl(1 +
\Omega^*\biggl(\sum_{n=1}^{\infty}\Lambda^n(V_{q-r})t^n\biggr)\biggr)
+ q\sum_{k=1}^\infty t^{kq}.$$
Thus, by \eqref{Hellerlog} and the definition of $\Lambda(V_{q-r},t)$,
\begin{eqnarray*}
\psi_S(V_r,t)\! &\ind& \! t\displaystyle\frac{d}{dt}\,\Omega^*
(\log \Lambda(V_{q-r},t)) + q\sum_{k=1}^\infty t^{kq}\\
&\ind& \!
\Omega^*\biggl(t\displaystyle\frac{d}{dt}\log \Lambda(V_{q-r},t)
\biggr) +q\sum_{k=1}^\infty t^{kq}.
\end{eqnarray*}
Therefore, by \eqref{AdamsdefLambda} and the definition of $\Omega^*$,
$$\psi_S(V_r,t) \ind
\sum_{n=1}^{\infty}(-1)^{n-1}\Omega^n(\psi_\Lambda^n(V_{q-r}))t^n +
q\sum_{k=1}^\infty t^{kq}.$$
Comparing coefficients of $t^n$ gives
\begin{equation}
\label{psiSind2} \psi_S^n(V_r) \ind \begin{cases}
(-1)^{n-1} \Omega^n(\psi_\Lambda^n(V_{q-r})) & \mbox{ if\, $q \nmid n$},\\
(-1)^{n-1} \Omega^n(\psi_\Lambda^n(V_{q-r})) + qV_1 & \mbox{ if\, $q
\mid n$},
\end{cases}
\end{equation}
for all $n>0$ and all $r$ such that $\tilde{q}\leq r\leq q$.

\begin{theorem}
\label{psiSLambda} Let $K$ be a field of prime characteristic $p$
and let $C$ be a cyclic $p$-group of order $q > 1$.
The Adams operations
$\psi_\Lambda^n$ and $\psi_S^n$ on $\Rkc$ satisfy
$$\psi_S^n(V_r) \proj (-1)^{n-1} \Omega^n(\psi_\Lambda^n(V_{q-r})) +
(n,q)V_{q/(n,q)}$$
for all $r$ such that\/ $q/p \leq r\leq q$ and all $n>0$.
\end{theorem}

Theorem \ref{psiSLambda} yields
$$\psi_S^n(V_r) = (-1)^{n-1} \Omega^n(\psi_\Lambda^n(V_{q-r})) +
(n,q)V_{q/(n,q)} + aV_q,$$
for some $a \in \mathbb{Z}$. If $\psi_\Lambda^n(V_{q-r})$ is known
then we may calculate $s=\delta(\Omega^n(\psi_\Lambda^n(V_{q-r})))$.
Since $\delta(\psi_S^n(V_r)) = r$, by \eqref{dimemsion}, we find that
$a = (r +(-1)^ns - q)/q$. Hence Theorem \ref{psiSLambda} determines
$\psi_S^n(V_r)$ completely in terms of $\psi_\Lambda^n(V_{q-r})$, for
$r \in \{q/p,\dots,q\}$. Suppose now that $r \in \{1,\dots,q\}$, and
choose $j$ such that $1 \leq p^{j-1} \leq r \leq p^j \leq q$. Then
$C$ has a factor group $C(p^j)$ of order $p^j$, and $R_{KC(p^j)}$
may be regarded as the subring of $\Rkc$ spanned by
$V_1,V_2,\dots,V_{p^j}$ (see the remarks preceding \eqref{Pmult}).
By Theorem \ref{psiSLambda} applied to $C(p^j)$ we may express
$\psi_S^n(V_r)$ in terms of $\psi_\Lambda^n(V_{p^j-r})$. Thus
Theorem \ref{psiSLambda} enables $\psi_S^n$ to be completely
determined from $\psi_\Lambda^n$.

\vspace{0.2in}
\noindent
{\bf Proof of Theorem \ref{psiSLambda}.}
As before, write $\tilde{q} = q/p$. We
prove the result by induction on $q$. First suppose that $q=p$ and
let $r \in \{1,\dots,p\}$. By \eqref{psiSind2} we have
$$\psi_S^n(V_r) \ind (-1)^{n-1} \Omega^n(\psi_\Lambda^n(V_{p-r})) +
(n,p)V_{p/(n,p)},$$
since $V_p$ is induced. The result then follows
from the fact that all induced $KC$-modules are projective, since $q=p$.

Now let $q>p$ and assume that the result holds for $\tilde{C}$.
Let $r \in \{\tilde{q},\dots,q\}$ and set
$$Y = \psi_S^n(V_r) + (-1)^n\Omega^n(\psi_\Lambda^n(V_{q-r})).$$
By \eqref{Adamsrestrict} and \eqref{Schanuel},
we have
$$Y\downa_{\tilde{C}} \ \proj \psi_S^n(V_r\downa_{\tilde{C}}) +
(-1)^n\Omega^n(\psi_\Lambda^n(V_{q-r}\downa_{\tilde{C}})).$$
Write $r = ap + b$ where $0 \leq b < p$. Then, by
\eqref{restriction} and \eqref{restrictheller}, we have
\begin{eqnarray*}
Y\downa_{\tilde{C}}\! &\proj& \!\!(p-b)\left[\psi_S^n(\tilde{V}_a) +
(-1)^n\Omega^n(\psi_\Lambda^n(\tilde{V}_{\tilde{q}-a}))\right]\\ &&
\qquad \quad \mbox{}+ b\left[\psi_S^n(\tilde{V}_{a+1}) +
(-1)^n\Omega^n(\psi_\Lambda^n(\tilde{V}_{\tilde{q}-(a+1)}))\right].
\end{eqnarray*}
Note that $\tilde{q}/p \leq a \leq \tilde{q}$ since $\tilde{q}
\leq r \leq q$. Also, if $b>0$ then $\tilde{q}/p <
a+1 \leq \tilde{q}$. Hence, by the inductive hypothesis, we may
replace each of the expressions contained in square brackets above
by $(n,\tilde{q})\tilde{V}_{\tilde{q}/(n,\tilde{q})}$.
Thus
\begin{equation}
\label{Yres}
Y\downa_{\tilde{C}} \ \proj
p(n,\tilde{q})\tilde{V}_{\tilde{q}/(n,\tilde{q})}.
\end{equation}

Consider the case where $q \nmid n$, so that $(n,q) = (n,\tilde{q}) < q$.
Then, by \eqref{psiSind2}, $Y$ is induced. Hence
$Y - (n,q)V_{q/(n,q)}$ is induced. Also, by \eqref{Yres},
$$(Y - (n,q)V_{q/(n,q)})\downa_{\tilde{C}} \ \proj
p(n,\tilde{q})\tilde{V}_{\tilde{q}/(n,\tilde{q})}
- p(n,q)\tilde{V}_{\tilde{q}/(n,q)} = 0.$$
Therefore $Y - (n,q)V_{q/(n,q)} \proj 0$ by Lemma \ref{resind}\,(ii).
This gives the required result.

Finally suppose that $q \mid n$, so that $(n,q) = q$ and
$(n,\tilde{q}) = \tilde{q}$. Then, by \eqref{psiSind2},
$Y - qV_1$ is induced. Also, by \eqref{Yres},
$$(Y - qV_1)\downa_{\tilde{C}} \ \proj p\tilde{q}\tilde{V}_1 -
q\tilde{V}_1 = 0.$$
Therefore $Y - qV_1 \proj 0$ by Lemma \ref{resind}\,(ii). This
gives the required result. \hfill $\square$

\vspace{0.2in}
As before $\delta$ denotes the endomorphism of $\Rkc$ satisfying
$\delta(V_r) = rV_1$ for all $r$.

\begin{corollary}
\label{periodL2} In the notation of Theorem\/ {\rm 6.1}, where $\sigma$
denotes the minimum period of\/ $\psi_S^n$, we have
$\psi_S^\sigma = \delta$.
\end{corollary}

\noindent
{\bf Proof.}
Since $\psi_S^\sigma = \psi_S^{2q}$, it is sufficient to show that
$\psi_S^{2q}(V_r) = rV_1$ for $r = 1,\dots,q$. Suppose first
that $r \in \{\tilde{q},\dots,q\}$. Then, by Theorem
\ref{psiSLambda}, we have
$$\psi_S^{2q}(V_r) \proj - \psi_\Lambda^{2q}(V_{q-r}) +qV_1.$$
Thus, by Theorem \ref{periodL}\,(iii), we have
$\psi_S^{2q}(V_r) \proj -(q-r)V_1 + qV_1 = rV_1$.
Since $\delta(\psi_S^{2q}(V_r)) = r$ it follows that
$\psi_S^{2q}(V_r) = rV_1$.

Now suppose that $r \leq \tilde{q}$.
Then we may regard $V_r$ as a module for the factor group
$C(\tilde{q})$ of $C$ with order $\tilde{q}$ and, by induction,
we may assume that $\psi_S^{2\tilde{q}}(V_r) = rV_1$. However, for
$R_{KC(\tilde{q})}$, the minimum period of $\psi_S^n$ divides
$2\tilde{q}$. Hence $\psi_S^{2q}(V_r) = rV_1$, as required.
\hfill $\square$

\section{Minimum periods}
\noindent
Let $G$ be a finite group and $K$ a field of characteristic
$p \geq 0$. Throughout this section we shall assume that the
Adams operations $\psi_\Lambda^n$ and $\psi_S^n$ are periodic in
$n$. Thus either $p=0$ or $p$ is a prime and the Sylow $p$-subgroups
of $G$ are cyclic. We shall give lower bounds for the minimum periods
of $\psi_\Lambda^n$ and $\psi_S^n$.

Let $G_{p'}$ denote the set of all elements of $G$ of order prime to
$p$, let $e$ denote the exponent of $G$,
and let $e'$ denote the least common multiple of the orders of the elements
of $G_{p'}$. Also, write $\lambda$ for the minimum period of
$\psi_\Lambda^n$ and $\sigma$ for the minimum period of $\psi_S^n$.
We begin by showing that $\lambda$ and $\sigma$ are divisible by
$e'$.

Let $V$ be the regular $KG$-module and let $\chi$ denote the
Brauer character of $V$. Thus, for all $g \in G_{p'}$, we have
$\chi(g) = |G|$ if $g=1$ and $\chi(g) = 0$ otherwise.
We extend the definition of Brauer character to elements of
$\Rkg$ by linearity (as in the proof of Lemma \ref{charzero})
and, for each $n > 0$, we write $\chi_\Lambda^n$
for the Brauer character of $\psi_\Lambda^n(V)$ and $\chi_S^n$ for
the Brauer character of $\psi_S^n(V)$. By \cite[Lemma 2.6]{Br2},
$\chi_\Lambda^n(g) = \chi(g^n) = \chi_S^n(g)$
for all $g \in G_{p'}$.  Since $\psi_\Lambda^{e'} =
\psi_\Lambda^{e'+\lambda}$ and $\psi_S^{e'} = \psi_S^{e'+\sigma}$
it follows that
$\chi(g^\lambda) = \chi(1) = \chi(g^\sigma)$
for all $g \in G_{p'}$. Thus $g^\lambda = g^\sigma = 1$ for all
$g \in G_{p'}$, giving that $e' \mid \lambda$ and $e' \mid \sigma$.

If $|G|$ is not divisible by $p$ then $G_{p'}= G$ and $e' = e$.
Thus, by the argument above
and Lemma \ref{charzero}, we have the following result.

\begin{proposition}
Let $G$ be a finite group and $K$ a field such that $|G|$ is not
divisible by the characteristic of\/ $K$. Then $\psi_\Lambda^n$ and
$\psi_S^n$ are both periodic in $n$ with minimum period $e$, where $e$
is the exponent of\/ $G$.
\end{proposition}

From now on we shall assume that $K$ has prime characteristic $p$.

\begin{proposition}
Let $K$ be a field of prime characteristic $p$ and let $G$ be a
finite group with a cyclic Sylow $p$-subgroup $C$ of order
$q > 1$. Then $\psi_\Lambda^n$ and $\psi_S^n$ are periodic in $n$
with minimum periods divisible by the least common multiple of\/ $2$
and the exponent of\/ $G$.
\end{proposition}

\noindent
{\bf Proof.}
Let $\lambda$, $\sigma$, $e$ and $e'$ be as before and let
$V$ be the regular $KG$-module. Then
$\psi_\Lambda^n(V\downa_C) = \psi_\Lambda^{n+\lambda}(V\downa_C)$ and
$\psi_S^n(V\downa_C) = \psi_S^{n+\sigma}(V\downa_C)$ for all
$n > 0$. Since $V\downa_C$ is a free $KC$-module, it follows
from Propositions 4.5 and 4.6 that $q \mid \lambda$ and
$q \mid \sigma$. But, as we have seen, $e' \mid \lambda$ and
$e' \mid \sigma$.
Hence $e \mid \lambda$ and $e \mid \sigma$. This
completes the proof for $p=2$. Thus we now assume that $p$ is odd.
Consider the $KC$-module $V_2$ and let $M = V_2\upa^G$. It is
straightforward to prove from Mackey's decomposition theorem (see
a similar argument in Section 3) that
$M\downa_C\ = rV_2 + W$, where $r$ is a positive integer
and $W$ is a sum of regular
modules for factor groups of $C$ regarded as $KC$-modules.
For all $n > 0$ we have
$$\psi_\Lambda^n(rV_2 + W) = \psi_\Lambda^{n+\lambda}(rV_2 + W),
\quad \psi_S^n(rV_2 + W) = \psi_S^{n+\sigma}(rV_2 + W).$$
However, by Lemma \ref{Adamsinflate}, $\psi_\Lambda^n(W) =
\psi_S^n(W) = W$ for all $n$ such that $p \nmid n$. Thus
$\psi_\Lambda^1(V_2) = \psi_\Lambda^{\lambda +1}(V_2)$ and
$\psi_S^1(V_2) = \psi_S^{\sigma +1}(V_2)$. By \cite[Theorem 5.1]{BJ}
it now follows that $\lambda$ and $\sigma$ are divisible by $2$.
\hfill $\square$

\vspace{0.2in}
The case $q=p$ was considered in \cite[Theorem 7.2]{Br2} and it
was shown that $\lambda \mid 2e$ and $\sigma \mid 2e$, where $e$
is the exponent of $G$.
Thus, in view of the results above and our results for the case $G=C$, we
make the following conjecture.

\vspace{0.1in}
\noindent
{\bf Conjecture.} Let $K$ be a field of prime characteristic $p$
and let $G$ be a finite group with cyclic Sylow $p$-subgroups.
Let $e$ be the exponent of\/ $G$. Then the minimum period of
$\psi_\Lambda^n$ is $e$ or $2e$, and the minimum period of
$\psi_S^n$ is $e$ or $2e$.

\vspace{0.1in}
\noindent
{\bf Acknowledgements.} Work supported by EPSRC Standard Research Grant
EP/G024898/1.
We are indebted to Peter Symonds for suggesting
the use of Conlon's induction theorem in Section 5.

\end{document}